\documentclass[journal]{IEEEtran}
\usepackage[utf8]{inputenc}
\usepackage[T1]{fontenc}

\usepackage{cite} 
\usepackage{graphicx}
\usepackage{subfigure}
\renewcommand{\thesubfigure}{\thefigure.\arabic{subfigure}}
\makeatletter
\renewcommand{\p@subfigure}{}
\renewcommand{\@thesubfigure}{\thesubfigure:\hskip\subfiglabelskip}
\makeatother

\usepackage{amsmath}
\usepackage{array}

\newcommand{\HH}{{\mathcal H}}
\newcommand{\II}{{\mathcal I}}
\newcommand{\KK}{{\mathcal K}}
\newcommand{\Sal}{{\mathcal S}}

\newcommand{\id}{i_d}

\newcommand{\OO}{{\mathcal O}}

\newcommand{\w}{{\omega}}
\newcommand{\W}{{\Omega}}
\newcommand{\eps}{{\varepsilon}}

\newcommand{\abs}[1]{\left\lvert#1\right\rvert}
\newcommand{\norm}[1]{\left\lVert#1\right\rVert}

\begin{document}
%
\title{Sensorless position estimation of Permanent-Magnet Synchronous Motors using a saturation model}
\author{Al~Kassem~Jebai, Fran\c{c}ois~Malrait, Philippe~Martin and~Pierre~Rouchon
\thanks{A-K.~Jebai, P.~Martin and P.~Rouchon are with the Centre Automatique et Syst\`{e}mes, MINES ParisTech, 75006~Paris,~France
{\tt\footnotesize \{al-kassem.jebai, philippe.martin, pierre.rouchon\}@mines-paristech.fr}}%
\thanks{{F.~Malrait is with Schneider Toshiba~Inverter~Europe}, 27120 Pacy-sur-Eure, France
{\tt\footnotesize francois.malrait@schneider-electric.com}}%
}

\maketitle

\begin{abstract}
Sensorless control of Permanent-Magnet Synchronous Motors (PMSM) at low velocity remains a challenging task. A now well-established method consists in injecting a high-frequency signal and use the rotor saliency, both geometric and magnetic-saturation induced. This paper proposes a clear and original analysis based on second-order averaging of how to recover the position information from signal injection; this analysis blends well with a general model of magnetic saturation. It also proposes a simple parametric model of the saturated PMSM, based on an energy function which simply encompasses saturation and cross-saturation effects. Experimental results on a surface-mounted PMSM and an interior magnet PMSM illustrate the relevance of the approach.
\end{abstract}


\section{Introduction}
\PARstart{P}{ermanent}-Magnet Synchronous Motors (PMSM) are widely used in industry. In the so-called ``sensorless'' mode of operation, the rotor position and velocity are not measured and the control law must make do with only current measurements. While sensorless control at medium to high velocities is well understood, with many reported control schemes and industrial products, sensorless control at low velocity remains a challenging task. The reason is that observability degenerates at zero velocity, causing a serious problem in the necessary rotor position estimation.

A now well-established method to overcome this problem is to add some persistent excitation by injecting a high-frequency signal~\cite{JanseL1995ITIA} and use the rotor saliency, whether geometric for Interior Permanent-Magnet machines or induced by main flux saturation for Surface Permanent-Magnet machines~\cite{OgasaA1998ITIA,CorleL1998ITIA,AiharTYME1999ITPE,ConsoST2001ITIA,HaISS2003ITIA,JangSHIS2003ITIA, JangHOIS2004ITIA,RobeiS2004ITIA,Shinn2008ITIA}.
Signal injection is moreover considered as a standard building block in hybrid control schemes for complete drives operating from zero to full speed~\cite{HarneN2000ITIE,WallmHC2005ITIE,SilvaAS2006ITIE,PiippHL2008ITIE,FooR2010ITIE}.

However to get a good position estimation under high-load condition it is important to take cross-saturation into account~\cite{GugliPV2006ITIA,BiancBJS2007ITPE,Holtz2008ITIA,ReigoGRBL2008ITIA,BiancBF2009IEMDC,DeKK2009ITPE,LiZHBS2009ITIA, SergeDM2009ITM,RacaGRBL2010ITIA,BiancFB2011ECCE,ZhuG2011ITIE}. It is thus necessary to rely on a model of the saturated PMSM adapted to control purposes, i.e. rich enough to capture in particular cross-saturation but also simple enough to be used in real-time and to be easily identified in the field; see \cite{ParasP1989ITEC,ChengCC2000ITM,StumbSDHT2003ITIA,RahmaH2005ITIA, StumbSSTD2005ITM,ArmanGPPV2009ITIA} for references more or less in this spirit.

The contribution of this paper, which builds on the preliminary work~\cite{JebaiMMR2011IEMDC}, is twofold: on the one hand it proposes a clear and original analysis based on second-order averaging of how to recover the position information from signal injection; this analysis can accommodate to any form of injected signals, e.g. square signals as in~\cite{YoonSMI2011ITIA}, and blends well with a general model of magnetic saturation including cross-saturation. On the other hand a simple parametric model of the saturated PMSM, well-adapted to control purposes, is introduced; it is based on an energy function which simply encompasses saturation and cross-saturation effects.

The paper runs as follows: section~\ref{sec:model} presents the saturation model. In section~\ref{sec:pos} position estimation by signal injection is studied thanks to second-order averaging. Section~\ref{sec:Identi} is devoted to the estimation of the parameters entering the saturation model using once again signal injection and averaging. Finally section~\ref{sec:experiment} experimentally demonstrates on two kinds of motors (with interior magnets and surface-mounted magnets) the relevance of the approach and the necessity of considering saturation to correctly estimate the position.

\section{An energy-based model of the saturated PMSM}\label{sec:model}

\subsection{Notations}
In the sequel we denote by $x_{ij}:=(x_i,x_j)^T$ the vector made from the real numbers $x_i$ and $x_j$, where
$ij$ can be $dq$, $\alpha\beta$ or~$\gamma\delta$. We also define the matrices
\[M_\mu:=\begin{pmatrix}\cos\mu& -\sin\mu\\ \sin\mu& \cos\mu\end{pmatrix}
\quad\text{and}\quad
\KK:=\begin{pmatrix}0&-1\\ 1&0\end{pmatrix}, \]
and we have the useful relation
\[\frac{dM_\mu}{d\mu}=\KK M_\mu=M_\mu\KK. \]

\subsection{Energy-based model}
The model of a two-axis PMSM expressed in the synchronous $d-q$ frame reads
\begin{align}
\label{eq:dqsys1}\frac{d\phi_{dq}}{dt} &=u_{dq}-Ri_{dq}-\w\KK(\phi_{dq}+\phi_{m})\\
\frac{J}{n^2}\frac{d\omega}{dt} &= \frac{3}{2}i_{dq}^T\KK(\phi_{dq}+\phi_{m}) - \frac{\tau_L}{n}\\
\label{eq:dqsys4}\frac{d\theta}{dt} &=\omega,
\end{align}
with $\phi_{dq}$ flux linkage due to the current; $\phi_{m}:=(\lambda,0)^T$ constant flux linkage due to the permanent magnet; $u_{dq}$ impressed voltage and $i_{dq}$ stator current; $\omega$ and $\theta$ rotor (electrical) speed and position; $R$ stator resistance; $n$ number of pole pairs; $J$ inertia moment and $\tau_L$ load torque. The physically impressed voltages are $u_{\alpha\beta}:=M_\theta u_{dq}$ while the physically measurable currents are~$i_{\alpha\beta}:=M_\theta i_{dq}$.
The current can be expressed in function of the flux linkage thanks to a suitable energy function~$\HH(\phi_d,\phi_q)$ by
\begin{align}
   \label{eq:CurrentFlux}i_{dq} =\II_{dq}(\phi_{dq})
   :=\begin{pmatrix}\partial_1\HH(\phi_d,\phi_q)\\\partial_2\HH(\phi_d,\phi_q)\end{pmatrix},
\end{align}
where $\partial_k\HH$ denotes the partial derivative w.r.t. the $k^{th}$ variable~\cite{BasicMR2010ITAC,BasicJMMR2011LNCIS}; without loss of generality $\HH(0,0)=0$. Such a relation between flux linkage and current naturally encompasses cross-saturation effects.

For an unsaturated PMSM this energy function reads
\begin{align*}
\HH_l(\phi_d,\phi_q) &=\frac{1}{2L_d}\phi^2_d+\frac{1}{2L_q}\phi^2_q
\end{align*}
where $L_d$ and $L_q$ are the motor self-inductances, and we recover the usual linear relations
\begin{align*}
    i_d &=\partial_1\HH_l(\phi_d,\phi_q) =\frac{\phi_d}{L_d}\\
    i_q &=\partial_2\HH_l(\phi_d,\phi_q) =\frac{\phi_q}{L_q}.
\end{align*}

Notice the expression for $\HH$ should respect the symmetry of the PMSM w.r.t the direct axis, i.e.
\begin{equation}\label{eq:sym}
    \HH(\phi_d,-\phi_q)=\HH(\phi_d,\phi_q),
\end{equation}
which is obviously the case for~$\HH_l$.
Indeed \eqref{eq:dqsys1}--\eqref{eq:dqsys4} is left unchanged by the transformation
\begin{multline*}
(u_d,u_q,\phi_d,\phi_q,i_d,i_q,\w,\theta,\tau_L)\rightarrow\\
(u_d,-u_q,\phi_d,-\phi_q,i_d,-i_q,-\w,-\theta,-\tau_L).
\end{multline*}

\subsection{Parametric description of magnetic saturation}\label{sec:paramsat}
Magnetic saturation can be accounted for by considering a more complicated magnetic energy function~$\HH$, having $\HH_l$ for quadratic part but including also higher-order terms. From experiments saturation effects are well captured by considering only third- and fourth-order terms, hence
\begin{multline*}
\HH(\phi_d,\phi_q)=\HH_l(\phi_d,\phi_q)\\
+\sum_{i=0}^3\alpha_{3-i,i}\phi_d^{3-i}\phi_q^i+\sum_{i=0}^4\alpha_{4-i,i}\phi_d^{4-i}\phi_q^i.
\end{multline*}
This is a perturbative model where the higher-order terms appear as corrections of the dominant term~$\HH_l$.
The nine coefficients $\alpha_{ij}$
together with $L_d$, $L_q$ are motor dependent. But \eqref{eq:sym} implies $\alpha_{2,1}=\alpha_{0,3}=\alpha_{3,1}=\alpha_{1,3}=0$, so that the energy function eventually reads
\begin{multline}\label{eq:EnerSat}
\HH(\phi_d,\phi_q) =\HH_l(\phi_d,\phi_q) +\alpha_{3,0}\phi_d^3+\alpha_{1,2}\phi_d\phi_q^2 \\
+\alpha_{4,0}\phi_d^4+\alpha_{2,2}\phi_d^2\phi_q^2+\alpha_{0,4}\phi_q^4.
\end{multline}
From~\eqref{eq:CurrentFlux} and~\eqref{eq:EnerSat} the currents are then explicitly given by
\begin{align}
 \label{eq:id} i_d &=\frac{\phi_d}{L_d}+3\alpha_{3,0}\phi_d^2+\alpha_{1,2}\phi_q^2 +4\alpha_{4,0}\phi_d^3+2\alpha_{2,2}\phi_d\phi_q^2\\
 \label{eq:iq} i_q &=\frac{\phi_q}{L_q}+2\alpha_{1,2}\phi_d\phi_q+2\alpha_{2,2}\phi_d^2\phi_q+4\alpha_{0,4}\phi_q^3,
\end{align}
which are the so-called flux-current magnetization curves.

To conclude, the model of the saturated PMSM is given by~\eqref{eq:dqsys1}--\eqref{eq:dqsys4} and~\eqref{eq:id}-\eqref{eq:iq}, with $\phi_d,\phi_q,\w,\theta$ as state variables. The magnetic saturation effects are represented by the five parameters $\alpha_{3,0},\alpha_{1,2},\alpha_{4,0},\alpha_{2,2},\alpha_{0,4}$.

\subsection{Model with $i_d,i_q$ as state variables}
The model of the PMSM is usually expressed with currents as state variables. This can be achieved here by time differentiating $i_{dq}=\II_{dq}(\phi_{dq})$,
\begin{align*}
    \frac{di_{dq}}{dt}
    =D\II_{dq}(\phi_{dq})\frac{d\phi_{dq}}{dt},
\end{align*}
with $\frac{d\phi_{dq}}{dt}$ given by~\eqref{eq:dqsys1}. Fluxes are then expressed as $\phi_{dq}=\II_{dq}^{-1}(i_{dq})$ by inverting the nonlinear relations~\eqref{eq:id}-\eqref{eq:iq}; rather than performing the exact inversion, we can take advantage of the fact the coefficients $\alpha_{i,j}$ are experimentally small. At first order w.r.t. the~$\alpha_{i,j}$ we have $\phi_d=L_di_d+\OO(\abs{\alpha_{i,j}})$ and $\phi_q=L_qi_q+\OO(\abs{\alpha_{i,j}})$; plugging these expressions into~\eqref{eq:id}-\eqref{eq:iq} and neglecting $\OO(\abs{\alpha_{i,j}}^2)$ terms, we easily find
\begin{align}
 \notag\phi_d &=L_d\bigl(i_d-3\alpha_{3,0}L_d^2i_d^2 -\alpha_{1,2}L_q^2i_q^2\\
 \label{eq:phid}&\quad\qquad-4\alpha_{4,0}L_d^3i_d^3-2\alpha_{2,2}L_dL_q^2i_di_q^2\bigr)\\
 \notag\phi_q &=L_q\bigl(i_q-2\alpha_{1,2}L_dL_qi_di_q-\\
 \label{eq:phiq}&\quad\qquad2\alpha_{2,2}L_d^2L_qi_d^2i_q-4\alpha_{0,4}L_q^3i_q^3\bigr).
\end{align}

Notice the matrix
\begin{equation}\label{eq:Gmat}
\begin{pmatrix}G_{dd}(i_{dq})&G_{dq}(i_{dq})\\G_{dq}(i_{dq})&G_{qq}(i_{dq})\end{pmatrix}
:=D\II_{dq}\bigl(\II_{dq}^{-1}(i_{dq})\bigr),
\end{equation}
with coefficients easily found to be
\begin{align*}
    G_{dd}(i_{dq})&=\frac{1}{L_d}+6\alpha_{3,0}L_di_d+12\alpha_{4,0}L_d^2i_d^2+2\alpha_{2,2}L_q^2i_q^2\\
    G_{dq}(i_{dq})&=2\alpha_{1,2}L_qi_q+4\alpha_{2,2}L_di_dL_qi_q\\
    G_{qq}(i_{dq})&=\frac{1}{L_q}+2\alpha_{1,2}L_di_d+2\alpha_{2,2}L_d^2i_d^2+12\alpha_{0,4}L_q^2i_q^2,
\end{align*}
is by construction symmetric; indeed
\begin{align*}
    D\mathcal{I}_{dq}(\phi_{dq})
    =\begin{pmatrix}\partial_{11}\HH(\phi_d,\phi_q) & \partial_{21}\HH(\phi_d,\phi_q)\\
    \partial_{12}\HH(\phi_d,\phi_q) & \partial_{22}\HH(\phi_d,\phi_q)\end{pmatrix}
\end{align*}
and $\partial_{12}\HH=\partial_{21}\HH$. Therefore the inductance matrix
\begin{align*}
   \begin{pmatrix}L_{dd}(i_{dq}) & L_{dq}(i_{dq})\\ L_{dq}(i_{dq})& L_{qq}(i_{dq})\end{pmatrix}
   :=&\begin{pmatrix}G_{dd}(i_{dq})&G_{dq}(i_{dq})\\G_{dq}(i_{dq})&G_{qq}(i_{dq})\end{pmatrix}^{-1}.
\end{align*}
is also symmetric, though this is not always acknowledged in saturation models encountered in the literature.

\section{Position estimation by high frequency voltage injection}\label{sec:pos}

\subsection{Signal injection and averaging}\label{sec:averaging}
A general sensorless control law can be expressed as
\begin{align}
    \label{eq:claw1}u_{\alpha\beta} &=M_{\theta_c}u_{\gamma\delta}\\
    \frac{d\theta_c}{dt} &=\omega_c\\
    \frac{d\eta}{dt} &=a\bigl(M_{\theta_c}i_{\gamma\delta},\theta_c,\eta,t\bigr)\\
    \label{eq:claw4}\omega_c &=\W_c\bigl(M_{\theta_c}i_{\gamma\delta},\theta_c,\eta,t\bigr)\\
    \label{eq:claw5}u_{\gamma\delta} &=\mathcal U_{\gamma\delta}\bigl(M_{\theta_c}i_{\gamma\delta},\theta_c,\eta,t\bigr),
\end{align}
where the measured currents $i_{\alpha\beta}=M_{\theta_c}i_{\gamma\delta}$ are used to compute $u_{\gamma\delta}$, $\omega_c$ and the evolution of the internal (vector) variable $\eta$ of the controller; $\theta_c$ and~$\omega_c$ are known by design.

It will be convenient to write the system equations~\eqref{eq:dqsys1}--\eqref{eq:dqsys4} in the $\gamma-\delta$ frame defined by $x_{\gamma\delta}:=M_{\theta-\theta_c}x_{dq}$, which gives
\begin{align}
\label{eq:gdsys1}\frac{d\phi_{\gamma\delta}}{dt}
&=u_{\gamma\delta}-Ri_{\gamma\delta}-\w_c\KK\phi_{\gamma\delta}
-\w\KK M_{\theta-\theta_c}\phi_{m}\\
\frac{J}{n^2}\frac{d\omega}{dt} &= \frac{3}{2}i_{\gamma\delta}^T\KK(\phi_{\gamma\delta}+M_{\theta-\theta_c}\phi_{m}) - \frac{\tau_L}{n}\\
\label{eq:gdsys4}\frac{d\theta}{dt} &=\omega,
\end{align}
where from~\eqref{eq:id}-\eqref{eq:iq} currents and fluxes are related by
\begin{equation}\label{eq:iphigd}
    i_{\gamma\delta}=M_{\theta-\theta_c}\mathcal{I}_{dq}(M^T_{\theta-\theta_c}\phi_{\gamma\delta}).
\end{equation}

To estimate the position we will superimpose on some desirable control law~\eqref{eq:claw5} a fast-varying pulsating voltage,
\begin{equation}\label{eq:HFvoltage}
    u_{\gamma\delta} =\mathcal U_{\gamma\delta}\bigl(M_{\theta_c}i_{\gamma\delta},\theta_c,\eta,t\bigr)
     + \widetilde u_{\gamma\delta}f(\Omega t),
\end{equation}
where $f$ is a $2\pi$-periodic function with zero mean and $\widetilde u_{\gamma\delta}$ could like~$\cal U_{\gamma\delta}$ depend on $M_{\theta_c}i_{\gamma\delta},\theta_c,\eta,t$ (though it is always taken constant in the sequel). The constant pulsation $\Omega$ is chosen ``large'', so that $f(\Omega t)$ can be seen as a ``fast'' oscillation; typically $\Omega:=2\pi\times500\,\text{rad/s}$ in the experiments in section~\ref{sec:experiment}.

If we apply this modified control law to~\eqref{eq:gdsys1}--\eqref{eq:gdsys4}, then it can be shown the solution of the closed loop system is
\begin{align}
    \label{eq:gdavrth1}\phi_{\gamma\delta} &= \overline\phi_{\gamma\delta} + \frac{\widetilde  u_{\gamma\delta}}{\Omega}F(\Omega t) + \OO{(\frac{1}{\Omega^2})}\\
    \omega &= \overline\omega + \OO{(\frac{1}{\Omega^2})} \\
    \theta &= \overline\theta + \OO{(\frac{1}{\Omega^2})} \\
    \theta_c &= \overline\theta_c + \OO{(\frac{1}{\Omega^2})} \\
    \label{eq:gdavrth5}\eta &=  \overline\eta + \OO{(\frac{1}{\Omega^2})},
\end{align}
where $F$ is the primitive of $f$ with zero mean ($F$ clearly has the same period as~$f$); $(\overline\phi_{\gamma\delta},\overline\omega,\overline\theta,\overline\theta_c,\overline\eta)$ is the ``slowly-varying'' component of $(\phi_{\gamma\delta},\omega,\theta,\theta_c,\eta)$, i.e. satisfies
\begin{align*}
\frac{d\overline\phi_{\gamma\delta}}{dt}
&=\overline u_{\gamma\delta}-R\overline i_{\gamma\delta}-\overline\w_c\KK\overline\phi_{\gamma\delta}
-\w\KK M_{\overline\theta-\overline\theta_c}\phi_{m}\\
\frac{J}{n^2}\frac{d\overline\omega}{dt} &=\frac{3}{2}\overline i_{\gamma\delta}^T\KK(\overline\phi_{\gamma\delta}+M_{\overline\theta-\overline\theta_c}\phi_{m})
-\frac{\tau_L}{n}\\
    \frac{d\overline\theta}{dt}& = \overline\omega\\
    \frac{d\overline\theta_c}{dt} &=\overline\omega_c\\
    \frac{d\overline\eta}{dt} &=a\bigl(M_{\overline\theta_c}\overline i_{\gamma\delta},\overline\theta_c,\overline\eta,t\bigr),
\end{align*}
where
\begin{align}
    \label{eq:bari}\overline i_{\gamma\delta} &=M_{\overline\theta-\overline\theta_c}\mathcal{I}_{dq}(M^T_{\overline\theta
    -\overline\theta_c}\overline\phi_{\gamma\delta})\\
    \notag\overline\omega_c &=\W_c\bigl(M_{\overline\theta_c}\overline i_{\gamma\delta},\overline\theta_c,\overline\eta,t\bigr)\\
    \notag\overline u_{\gamma\delta} &=\mathcal U_{\gamma\delta}\bigl(M_{\overline\theta_c}\overline i_{\gamma\delta},\overline\theta_c,\overline\eta,t\bigr).
\end{align}
Notice this slowly-varying system is exactly the same as~\eqref{eq:gdsys1}--\eqref{eq:gdsys4} acted upon by the unmodified control law~\eqref{eq:claw1}--\eqref{eq:claw5}. In other words adding signal injection:
\begin{itemize}
  \item has a very small effect of order $\OO{(\frac{1}{\Omega^2})}$ on the mechanical variables $\theta,\omega$ and the controller variables $\theta_c,\eta$
  \item has a small effect of order $\OO{(\frac{1}{\Omega})}$ on the flux~$\phi_{\gamma\delta}$; this effect will be used in the next section to extract the position information from the measured currents.
\end{itemize}

The proof relies on a direct application of second-order averaging of differential equations, see~\cite{SandeVM2007book} section~$2.9.1$ and for the slow-time dependance section~$3.3$. Indeed setting $\eps:=\frac{1}{\W}$, $\sigma:=\frac{t}{\eps}$, and $x:=(\phi_{\gamma\delta},\omega,\theta,\theta_c,\eta)$, \eqref{eq:gdsys1}--\eqref{eq:gdsys4} acted upon by the modified control law \eqref{eq:claw1}--\eqref{eq:claw4} and~\eqref{eq:HFvoltage} is in the so-called standard form for averaging (with slow-time dependance)
\begin{align*}
    \frac{dx}{d\sigma} &=\eps f_1(x,\eps\sigma,\sigma)
    :=\eps\bigl(\overline f_1(x,\eps\sigma)+\widetilde f_1(x,\eps\sigma)f(\sigma)\bigr),
\end{align*}
with $f_1$ $T$-periodic w.r.t. its third variable ($T=2\pi$ in our case) and $\eps$ as a small parameter.  Therefore its solution can be approximated as 
\begin{align*}
    x(\sigma)&=z(\sigma)+\eps\bigl(u_1(z(\sigma),\eps\sigma,\sigma\bigr)+\OO(\eps^2),
\end{align*}
where $z(\sigma)$ is the solution of
\begin{align*}
    \frac{dz}{d\sigma} &=\eps g_1(z,\eps\sigma) + \eps^2g_2(z,\eps\sigma)
\end{align*}
and
\begin{IEEEeqnarray*}{rCl}
    g_1(y,\eps\sigma) &:=& \frac{1}{T}\int_0^Tf_1(y,\eps\sigma,s)ds
    =\overline f_1(y,\eps\sigma)\\
    v_1(y,\eps\sigma,\sigma) &:=& \int_0^\sigma\bigl(f_1(y,\eps\sigma,s)-g_1(y,\eps\sigma)\bigr)ds\\
    &=& \widetilde f_1(y,\eps\sigma)\int_0^\sigma f(s)ds\\
    u_1(y,\eps\sigma,\sigma) &:=& v_1(y,\eps\sigma,\sigma)-\frac{1}{T}\int_0^Tv_1(y,\eps\sigma,s)ds\\
    &=& \widetilde f_1(y,\eps\sigma)F(\sigma)\\
    K_2(y,\eps\sigma,\sigma) &:=& \partial_1f_1(y,\eps\sigma,\sigma)u_1(y,\eps\sigma,\sigma)\\
    &&{}-\partial_1u_1(y,\eps\sigma,\sigma)g_1(y,\eps\sigma)\\
    &=& [\overline f_1,\widetilde f_1](y,\eps\sigma)F(\sigma)\\
    &&{}+\frac{1}{2}\partial_1\widetilde f_1(y,\eps\sigma)\widetilde f_1(y,\eps\sigma)
    \frac{dF^2(\sigma)}{d\sigma}\\
    g_2(y,\eps\sigma) &:=& \frac{1}{T}\int_0^TK_2(y,\eps\sigma,s)ds=0.
\end{IEEEeqnarray*}
We have set
\[[\overline f_1,\widetilde f_1](y,\eps\sigma)\!:=\!\partial_1\overline f_1(y,\eps\sigma)\widetilde f_1(y,\eps\sigma)-\partial_1\widetilde f_1(y,\eps\sigma)\overline f_1(y,\eps\sigma)\]
and $F(\sigma):=\int_0^\sigma f(s)ds-\frac{1}{T}\int_0^T\int_0^\sigma f(s)dsd\sigma$, i.e. $F$ is the (of course $T$-periodic) primitive of~$f$ with zero mean.

Translating back to the original variables eventually yields the desired result~\eqref{eq:gdavrth1}--\eqref{eq:gdavrth5}.

\subsection{Position estimation}\label{sec:posesti}
We now express the effect of signal injection on the currents: plugging~\eqref{eq:gdavrth1} into~\eqref{eq:iphigd} we have
\begin{align}
    \notag i_{\gamma\delta}
    &=M_{\overline\theta-\overline\theta_c+\OO{(\frac{1}{\Omega^2})}}\\
    \notag&\qquad\mathcal{I}_{dq}\Bigl(M^T_{\overline\theta-\overline\theta_c+\OO{(\frac{1}{\Omega^2})}}
    \bigl(\overline\phi_{\gamma\delta}
    + \frac{\widetilde  u_{\gamma\delta}}{\Omega}F(\Omega t) + \OO{(\frac{1}{\Omega^2})}\bigr)\Bigr)\\
    \label{eq:ilowhigh}&=\overline i_{\gamma\delta}+\widetilde i_{\gamma\delta}F(\Omega t)+\OO{(\frac{1}{\Omega^2})},
\end{align}
where we have used~\eqref{eq:bari} and performed a first-order expansion to get
\begin{align}
    \notag\widetilde i_{\gamma\delta}
    &:=M_{\overline\theta-\overline\theta_c}D\mathcal{I}_{dq}\bigl(M^T_{\overline\theta-\overline\theta_c}
    \overline\phi_{\gamma\delta}\bigr)M^T_{\overline\theta-\overline\theta_c}\frac{\widetilde  u_{\gamma\delta}}{\Omega}\\
    \label{eq:tildei}&=M_{\overline\theta-\overline\theta_c}D\mathcal{I}_{dq}
    \Bigl(\mathcal{I}_{dq}^{-1}\bigl(M^T_{\overline\theta-\overline\theta_c}
    \overline i_{\gamma\delta}\bigr)\Bigr)M^T_{\overline\theta-\overline\theta_c}\frac{\widetilde  u_{\gamma\delta}}{\Omega}.
\end{align}
We will see in the next section how to recover $\widetilde i_{\gamma\delta}$ and $\overline i_{\gamma\delta}$ from the measured currents~$i_{\gamma\delta}$.
Therefore \eqref{eq:tildei} gives two (redundant) relations relating the unknown angle~$\overline\theta$ to the known variables $\overline\theta_c,\widetilde i_{dq},\overline i_{\gamma\delta},\widetilde u_{dq}$, provided the matrix
\[ \Sal(\mu,\overline i_{\gamma\delta}):=M_\mu D\mathcal{I}_{dq}
    \Bigl(\mathcal{I}_{dq}^{-1}\bigl(M^T_\mu\overline i_{\gamma\delta}\bigr)\Bigr)M^T_\mu \]
effectively depends on its first argument~$\mu$. This ``saliency condition'' is what is needed to ensure nonlinear observability.
The explicit expression for $\Sal(\mu,\overline i_{\gamma\delta})$ is obtained thanks to~\eqref{eq:Gmat}. In the case of an unsaturated magnetic circuit this matrix boils down to
\begin{align*}
    \Sal(\mu,\overline i_{\gamma\delta}) &=M_\mu
    \begin{pmatrix}\frac{1}{L_d}&0\\0&\frac{1}{L_q}\end{pmatrix}M^T_\mu\\
    &=\textstyle\frac{L_d+L_q}{2L_dL_q}
    \begin{pmatrix}1+\frac{L_d-L_q}{L_d+L_q}\cos2\mu &\frac{L_d-L_q}{L_d+L_q}\sin2\mu\\
    \frac{L_d-L_q}{L_d+L_q}\sin2\mu &1-\frac{L_d-L_q}{L_d+L_q}\cos2\mu\end{pmatrix}
\end{align*}
and does not depend on~$i_{\gamma\delta}$; notice this matrix does not depend on~$\mu$ for an unsaturated machine with no geometric saliency. Notice also \eqref{eq:tildei} defines in that case two solutions on~$]-\pi,\pi]$ for the angle~$\overline\theta$ since $\Sal(\mu,\overline i_{\gamma\delta})$ is actually a function of~$2\mu$; in the saturated case there is generically only one solution, except for some particular values of~$\overline i_{\gamma\delta}$.

There are several ways to extract the rotor angle information from~\eqref{eq:tildei}, especially for real-time use inside a feedback law. In this paper we just want to demonstrate the validity of~\eqref{eq:tildei} and we will be content with directly solving it through a nonlinear least square problem; in other words we estimate the rotor position as
\begin{equation}\label{eq:nonlin}
      \widehat\theta = \theta_c + \arg\min_{\mu\in]-\pi,\pi]}\norm{
      \widetilde i_{\gamma\delta} - \Sal(\mu,\overline i_{\gamma\delta})\frac{\widetilde u_{\gamma\delta}}{\Omega}}^2.
\end{equation}

\subsection{Current demodulation}\label{sec:currdem}
To estimate the position information using e.g.~\eqref{eq:nonlin} it is necessary to extract the low- and high-frequency components $\overline i_{\gamma\delta}$ and $\widetilde i_{\gamma\delta}$ from the measured current~$i_{\gamma\delta}$.
Since by~\eqref{eq:ilowhigh} $i_{\gamma\delta}(t)\approx\overline i_{\gamma\delta}(t)+\widetilde i_{\gamma\delta}(t)F(\Omega t)$ with $\overline i_{\gamma\delta}$ and $\widetilde i_{\gamma\delta}$ by construction nearly constant on one period of~$F$, we may write
\begin{align*}
    \overline i_{\gamma\delta}(t) &\approx\frac{1}{T}\int_{t-T}^ti_{\gamma\delta}(s)ds\\
    \widetilde i_{\gamma\delta}(t) &
    \approx\frac{\int_{t-T}^ti_{\gamma\delta}(s)F(\Omega s)ds}{\int_0^TF^2(\Omega s)ds},
\end{align*}
where $T:=\frac{2\pi}{\Omega}$. Indeed as $F$ is $2\pi$-periodic with zero mean,
\begin{align*}
    \int_{t-T}^ti_{\gamma\delta}(s)ds &\approx
    \overline i_{\gamma\delta}(t)\int_{t-T}^tds + \widetilde i_{\gamma\delta}(t)\int_{t-T}^tF(\Omega s)ds\\
    &=T\overline i_{\gamma\delta}(t)\\
    \int_{t-T}^ti_{\gamma\delta}(s)F(\Omega s)ds &\approx
    \overline i_{\gamma\delta}(t)\int_{t-T}^tF(\Omega s)ds \\
    &\qquad+ \widetilde i_{\gamma\delta}(t)\int_{t-T}^tF^2(\Omega s)ds\\
    &=\widetilde i_{\gamma\delta}(t)\int_0^TF^2(\Omega s)ds.
\end{align*}

\begin{figure}[ht!]
\centering
\includegraphics[width=0.9\columnwidth]{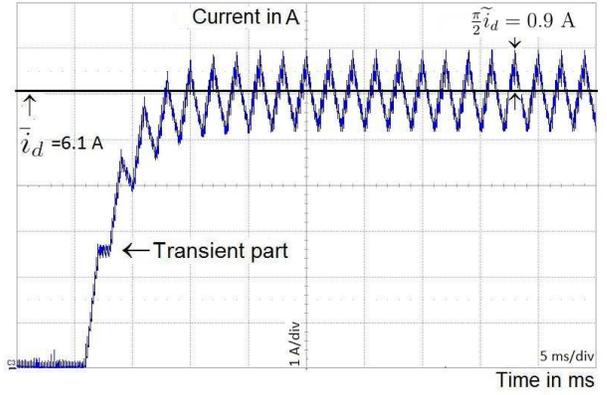}
\caption{Experimental time response of~$i_d$ in~\eqref{eq:dqlocked}-\eqref{eq:dqIF}}
\label{fig:HFcurrent}
\end{figure}
\section{Estimation of magnetic parameters}\label{sec:Identi}
The seven parameters in the saturation model~\eqref{eq:id}-\eqref{eq:iq} must of course be estimated. This can be done with a rather simple procedure also relying on signal injection and averaging.

\subsection{Principle}
The rotor is locked in the position~$\theta:=0$, hence the model \eqref{eq:dqsys1}--\eqref{eq:dqsys4} reduces to $\omega=0$ and
\begin{IEEEeqnarray}{rCl}
\frac{d\phi_{dq}}{dt} &=& u_{dq}-Ri_{dq}\label{eq:dqlocked},
\end{IEEEeqnarray}
with~$i_{dq}=\mathcal{I}_{dq}(\phi_{dq})$. Moreover $u_{dq}$ can now be physically impressed and $i_{dq}$ physically measured.

As in section~\ref{sec:averaging}, but now working directly in the $d-q$ frame, we inject a fast-varying pulsating voltage \begin{IEEEeqnarray}{rCl}
u_{dq} &=& \overline u_{dq}+\widetilde u_{dq}f(\W t)\label{eq:dqIF},
\end{IEEEeqnarray}
with constant $\overline u_{dq}$ and~$\widetilde u_{dq}$. The solution of~\eqref{eq:dqlocked}-\eqref{eq:dqIF} is then
\begin{IEEEeqnarray*}{rCl}
\phi_{dq} &=& \overline\phi_{dq} + \frac{\widetilde  u_{dq}}{\Omega}F(\Omega t) + \OO(\frac{1}{\Omega^2})
\end{IEEEeqnarray*}
where $\overline\phi_{dq}$, the ``slowly-varying'' component of $\phi_{dq}$, satisfies
\begin{IEEEeqnarray}{rCl}
\frac{d\overline\phi_{dq}}{dt} & = & \overline u_d-R\overline i_d\label{eq:dqslow},
\end{IEEEeqnarray}
with $\overline i_{dq}=\mathcal{I}_{dq}(\overline\phi_{dq})$. Moreover \eqref{eq:tildei} now boils down to
\begin{IEEEeqnarray}{rCl}
    \widetilde i_{dq} &=& D\mathcal{I}_{dq}\bigl(\mathcal{I}_{dq}^{-1}(\overline i_{dq})\bigr)\frac{\widetilde  u_{dq}}{\Omega}\label{eq:dqest}.
\end{IEEEeqnarray}

Since $\overline u_{dq}$ is constant \eqref{eq:dqslow} implies $R\overline i_{dq}$ tends to $\overline u_{dq}$, hence after an initial transient $\overline i_{dq}$ is constant. As a consequence $\widetilde i_{dq}$ is by~\eqref{eq:dqest} also constant. Fig.~\ref{fig:HFcurrent} shows for instance the time response of~$i_d$ for the SPM motor of section~\ref{sec:experiment} starting from~$\id(0)=0$ and using a square function~$f$; notice the current ripples seen on the scope are
$\max_{\tau\in[0,2\pi]}F(\tau)=\frac{\pi}{2}$ (since $f$ is square with period~$2\pi$) smaller than~$\widetilde i_{dq}$.

The magnetic parameters can then be estimated repeatedly using~\eqref{eq:dqest} with various values of~$\overline u_{dq}$ and $\widetilde u_{dq}$, as detailed in the next section.

\begin{table}[t!]
\setlength{\extrarowheight}{2pt}
\caption{Rated and estimated magnetic parameters of test motors\label{tab:param}}
\centering
\begin{tabular}{| l | l | l | }
  \firsthline
  Motor & IPM & SPM\\ \hline
  Rated power & 750~W & 1500~W\\ \hline
  Rated current $I_n$ (peak) & 4.51~A & 5.19~A\\ \hline
  Rated voltage (peak per phase) & 110~V & 245~V\\ \hline
  Rated speed & 1800~rpm & 3000~rpm\\ \hline
  Rated torque & 3.98~Nm & 6.06~Nm\\ \hline
  $n$ & 3 & 5\\ \hline
  $R$ & $1.52~\W$ & $2.1~\W$ \\ \hline
  $\lambda$ (peak) & $196$~mWb & $155$~mWb \\  \hline\hline
  $L_d$& $9.15$~mH &  $7.86$~mH\\ \hline
  $L_q$& $13.58$~mH & $8.18$~mH\\ \hline
  $\alpha_{3,0}L_d^2I_n$& $0.039$ & $0.056$\\ \hline
  $\alpha_{1,2}L_dL_qI_n$& $0.053$ & $0.055$\\ \hline
  $\alpha_{4,0}L_d^3I_n^2$& $0.0051$ & $0.0164$\\ \hline
  $\alpha_{2,2}L_dL_q^2I_n^2$& $0.0171$ & $0.027$\\ \hline
  $\alpha_{0,4}L_q^3I_n^2$& $0.0060$ & $0.0067$\\
  \lasthline
\end{tabular}
\end{table}
\begin{figure}[b!]
\centering
\includegraphics[width=0.9\columnwidth]{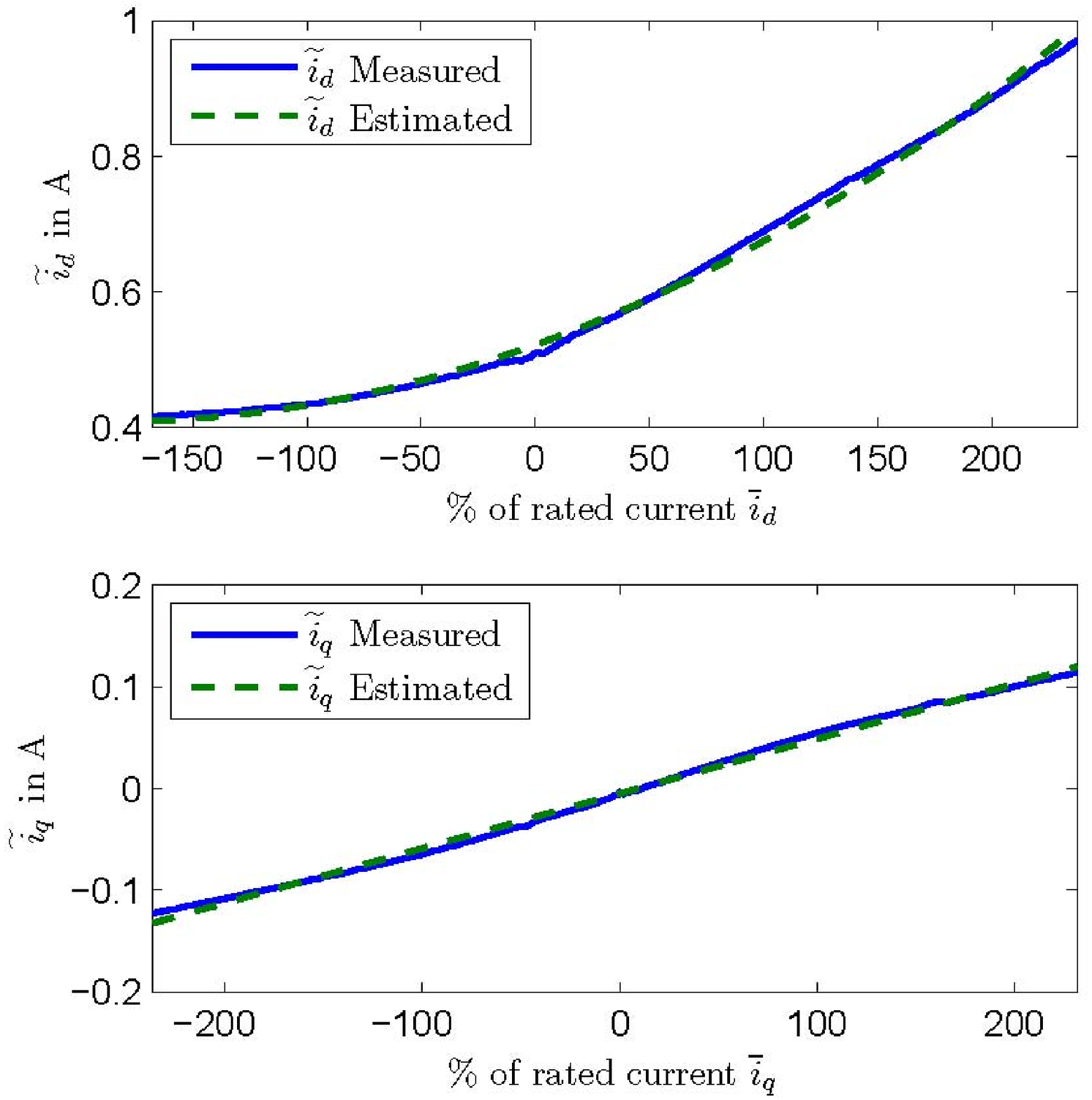}
\caption{IPM: fitted values vs measurements for \eqref{eq:saturation1} and~\eqref{eq:saturation3q}}
\label{fig:i_udIPM}
\end{figure}
\begin{figure}[b!]
\centering
\includegraphics[width=0.9\columnwidth]{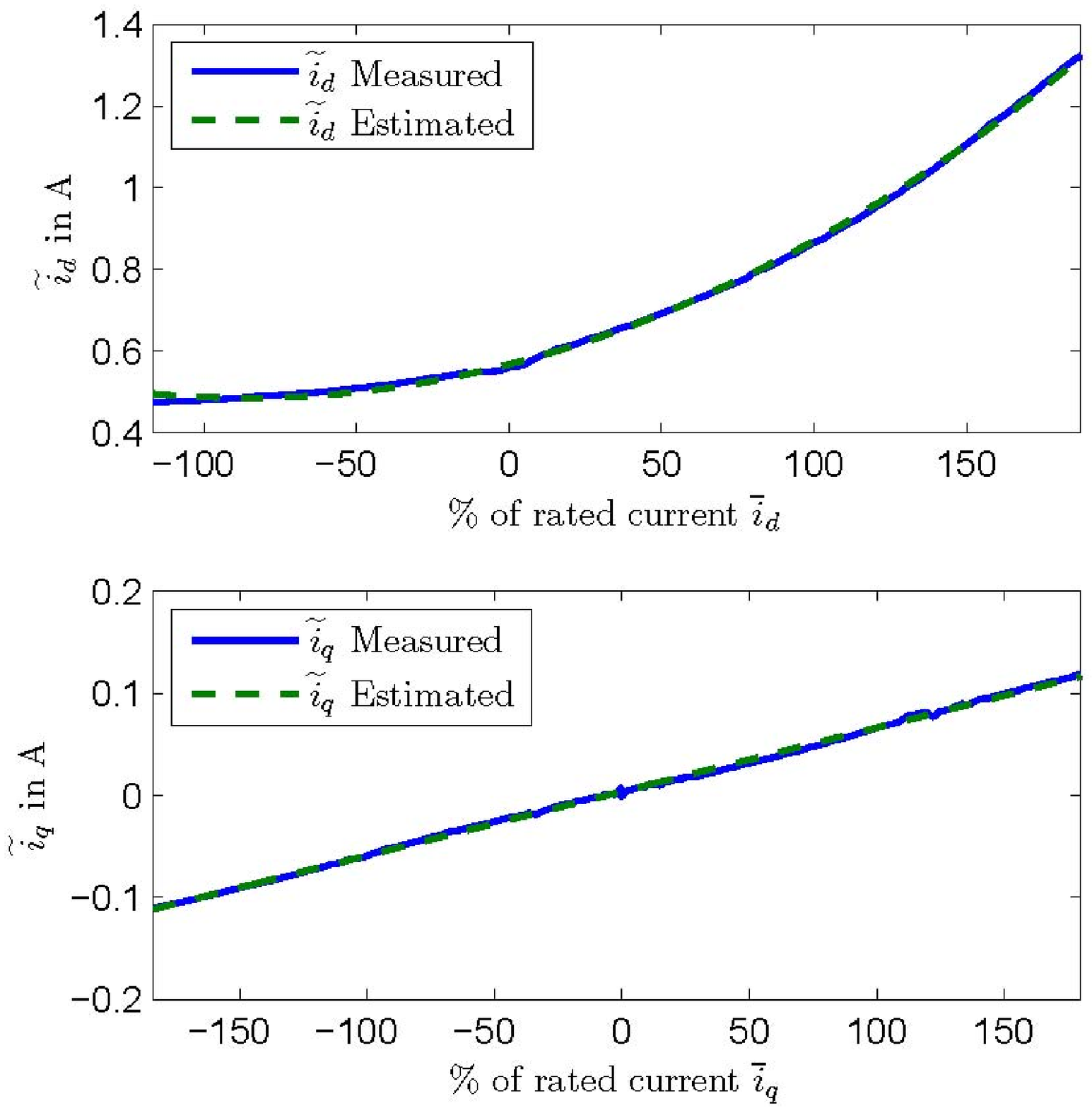}
\caption{SPM: fitted values vs measurements for \eqref{eq:saturation1} and~\eqref{eq:saturation3q}}
\label{fig:i_udSPM}
\end{figure}
\subsection{Estimation of the parameters}\label{sec:estpar}
From~\eqref{eq:Gmat} the entries of~$D\mathcal{I}_{dq}\bigl(\mathcal{I}_{dq}^{-1}(\overline i_{dq})\bigr)$ are given by
\begin{align*}
    G_{dd}(\overline i_{dq})&=\frac{1}{L_d}+6\alpha_{3,0}L_d\overline i_d+12\alpha_{4,0}L_d^2\overline i_d^2+2\alpha_{2,2}L_q^2\overline i_q^2\\
    G_{dq}(\overline i_{dq})&=2\alpha_{1,2}L_q\overline i_q+4\alpha_{2,2}L_d\overline i_dL_q\overline i_q\\
    G_{qq}(\overline i_{dq})&=\frac{1}{L_q}+2\alpha_{1,2}L_d\overline i_d+2\alpha_{2,2}L_d^2\overline i_d^2+12\alpha_{0,4}L_q^2\overline i_q^2.
\end{align*}
Since combinations of the magnetic parameters always enter linearly those equations, they can be estimated by simple linear least squares; moreover by suitably choosing $\overline u_{dq}$ and $\widetilde u_{dq}$, the whole least squares problem for the seven parameters can be split into several subproblems involving fewer parameters:
\begin{itemize}
 \item with $\overline u_{dq}:=0$, hence $\overline i_{dq}=0$, \eqref{eq:dqest} reads
\begin{IEEEeqnarray}{rCl}
    L_d &=& \frac{1}{\W}\frac{\widetilde u_d}{\widetilde i_d}\label{eq:Ld}\\
    L_q &=& \frac{1}{\W}\frac{\widetilde u_q}{\widetilde i_q}\label{eq:Lq}
\end{IEEEeqnarray}
 \item with $\overline u_q=0$, hence $\overline i_q=0$, and $\widetilde u_q=0$ \eqref{eq:dqest} reads
\begin{IEEEeqnarray}{rCl}
    \widetilde i_d &=& \frac{\widetilde u_d}{\W}\left(\frac{1}{L_d}+6\alpha_{3,0}L_d\overline i_d+12\alpha_{4,0}L_d^2\overline i_d^2\right)\label{eq:saturation1}\\
    \widetilde i_q &=& 0\notag
\end{IEEEeqnarray}
 \item with $\overline u_d=0$, hence $\overline i_d:=0$, and $\widetilde u_q=0$ \eqref{eq:dqest} reads
\begin{align}
   \label{eq:saturation3d}\widetilde i_d  &=\frac{\widetilde u_d}{\W}\Bigl(\frac{1}{L_d}+2\alpha_{2,2}L_q^2\overline i_q^2\Bigr)\\
   \label{eq:saturation3q}\widetilde i_q & =\frac{2\widetilde u_d}{\W}\alpha_{1,2}L_q\overline i_q
\end{align}
 \item with $\overline u_d=0$, hence $\overline i_d:=0$, and $\widetilde u_d=0$ \eqref{eq:dqest} reads
\begin{align}
  \label{eq:saturation4d}\widetilde i_d &=\frac{2\widetilde u_q}{\W}\alpha_{1,2}L_q\overline i_q\\
  \label{eq:saturation4q}\widetilde i_q &=\frac{\widetilde u_q}{\W}\Bigl(\frac{1}{L_q}+12\alpha_{0,4}L_q^2\overline i_q^2\Bigr).
\end{align}
\end{itemize}

$L_d$ and~$L_q$ are then immediately determined from \eqref{eq:Ld} and~\eqref{eq:Lq}; $\alpha_{3,0}$ and $\alpha_{4,0}$ are jointly estimated by least squares from~\eqref{eq:saturation1}; $\alpha_{2,2}$, $\alpha_{1,2}$ and $\alpha_{0,4}$ are separately estimated by least squares from respectively \eqref{eq:saturation3d}, \eqref{eq:saturation3q}-\eqref{eq:saturation4d} and~\eqref{eq:saturation4q}.

\begin{figure}[b!]
\centering
\includegraphics[width=0.9\columnwidth]{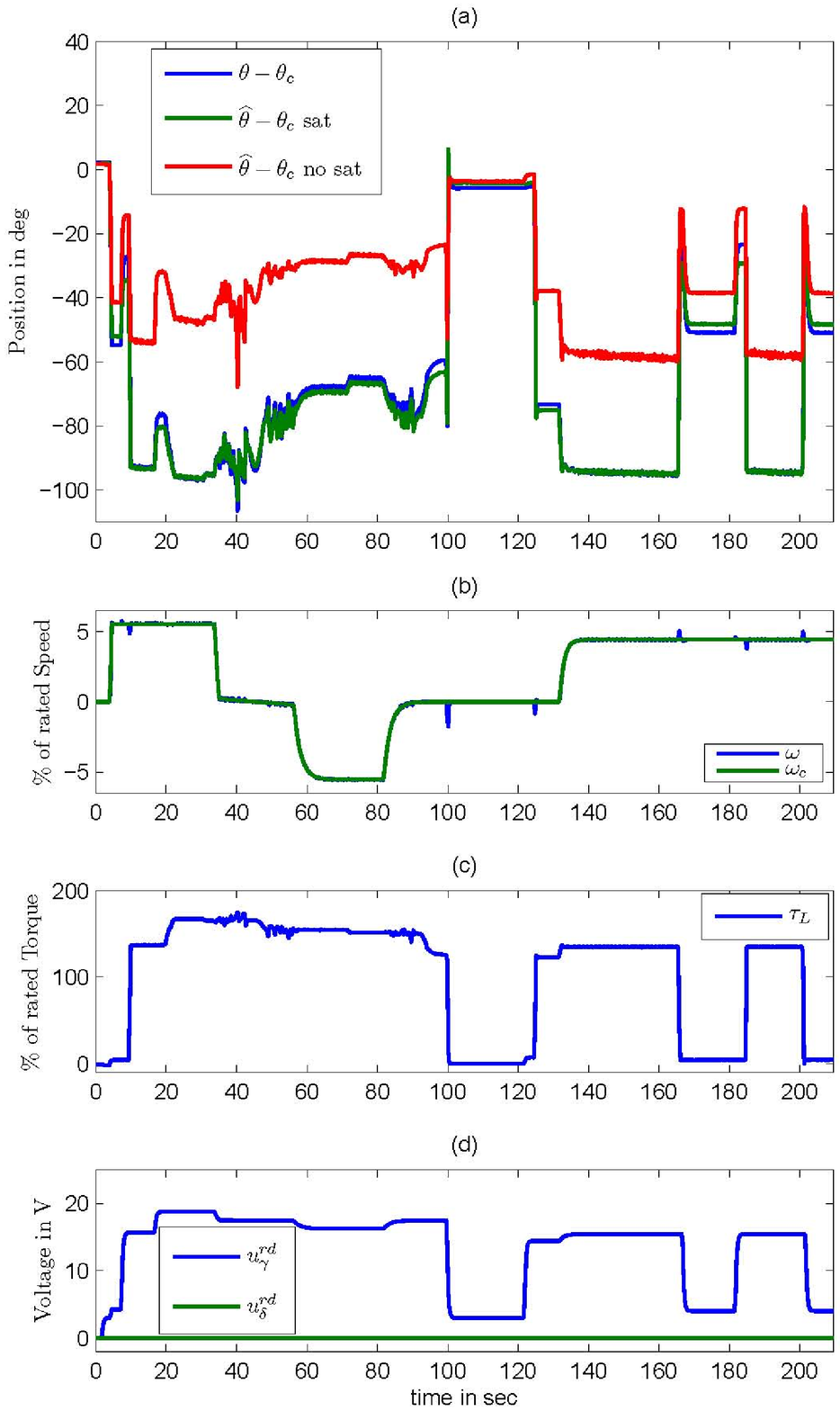}
\caption{Long test under various conditions for IPM: (a) measured $\theta-\theta_c$, estimated $\widehat\theta-\theta_c$ with and without saturation model; (b) measured speed $\omega$, reference speed~$\omega_c$; (c) load torque $\tau_L$; (d) voltages $u_{\gamma\delta}^{rd}$}
\label{fig:LongTestIPM}
\end{figure}
\begin{figure}[b!]
\centering
\includegraphics[width=0.9\columnwidth]{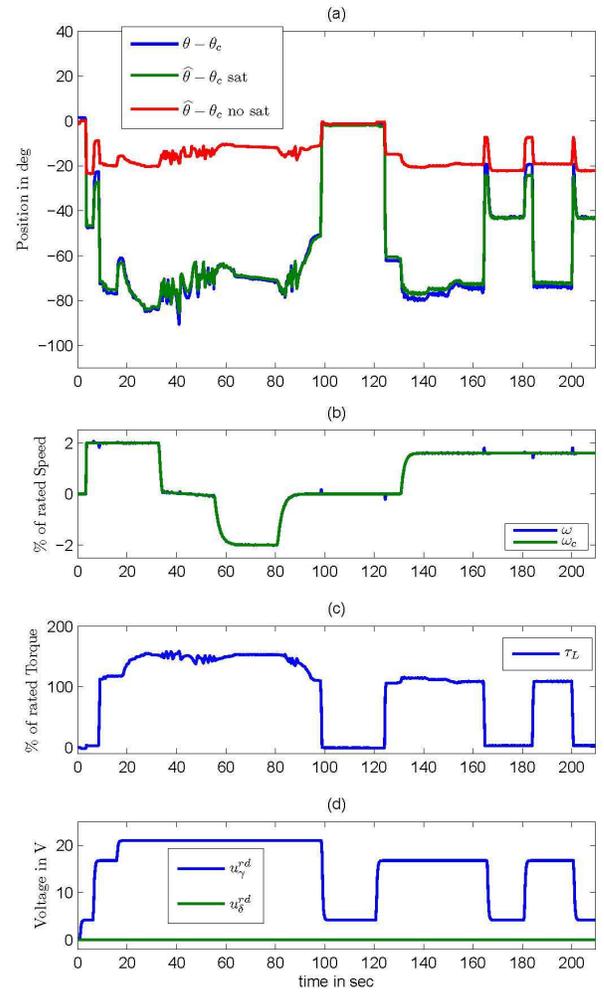}
\caption{Long test under various conditions for SPM: (a) measured $\theta-\theta_c$, estimated $\widehat\theta-\theta_c$ with and without saturation model; (b) measured speed $\omega$, reference speed~$\omega_c$; (c) load torque $\tau_L$; (d) voltages $u_{\gamma\delta}^{rd}$}
\label{fig:LongTestSPM}
\end{figure}
\subsection{Experimental setup} \label{sec:experiment}
The methodology developed in the paper was tested on two types of motors, an Interior Magnet PMSM (IPM) and a Surface-Mounted PMSM (SPM), with rated parameters listed in the top part of table~\ref{tab:param}.

The experimental setup consists of an industrial inverter ($400$\,V DC bus, $4$\,kHz PWM frequency), an incremental enco\-der, a dSpace fast prototyping system with 3 boards (DS1005, DS5202 and EV1048), and a host~PC. The measurements are sampled also at~$4$\,kHz, and synchronized with the PWM frequency. The load torque is created by a $4$\,kW DC motor.

\subsection{Estimation of the magnetic parameters}
We follow the procedure described in section~\ref{sec:Identi}: with the rotor locked in the position~$\theta:=0$, a square wave voltage with frequency $\Omega:=2\pi\times500\,\text{rad/s}$ and constant amplitude $\widetilde u_d$ or $\widetilde u_q$ (15\,V for the IPM, $14$\,V for the SPM) is applied to the motor; but for the determination of $L_d,L_q$ where $\overline u_d=\overline u_q:=0$, several runs are performed with various $\overline u_d$ (resp.~$\overline u_q$) such that $\overline i_d$ (resp.~$\overline i_q$) ranges from $-200\%$ to $+200\%$ of the rated current. The magnetic parameters are then estimated by linear least squares according to section~\ref{sec:estpar}, yielding the values in the bottom part of table~\ref{tab:param}.
Notice the SPM exhibits as expected little geometric saliency ($L_d\approx L_q$) hence the saturation-induced saliency is paramount to estimate the rotor position. Notice also the cross-saturation term~$\alpha_{12}$ is as expected quantitatively important for both motors.

The good agreement between the fitted curves and the measurements is demonstrated for instance for \eqref{eq:saturation1} and~\eqref{eq:saturation3q} on Fig.~\ref{fig:i_udIPM}-\ref{fig:i_udSPM}; notice \eqref{eq:saturation1} corresponds to saturation on a single axis while \eqref{eq:saturation3q} corresponds to cross-saturation.

\begin{figure}[b!]
\centering
\includegraphics[width=0.9\columnwidth]{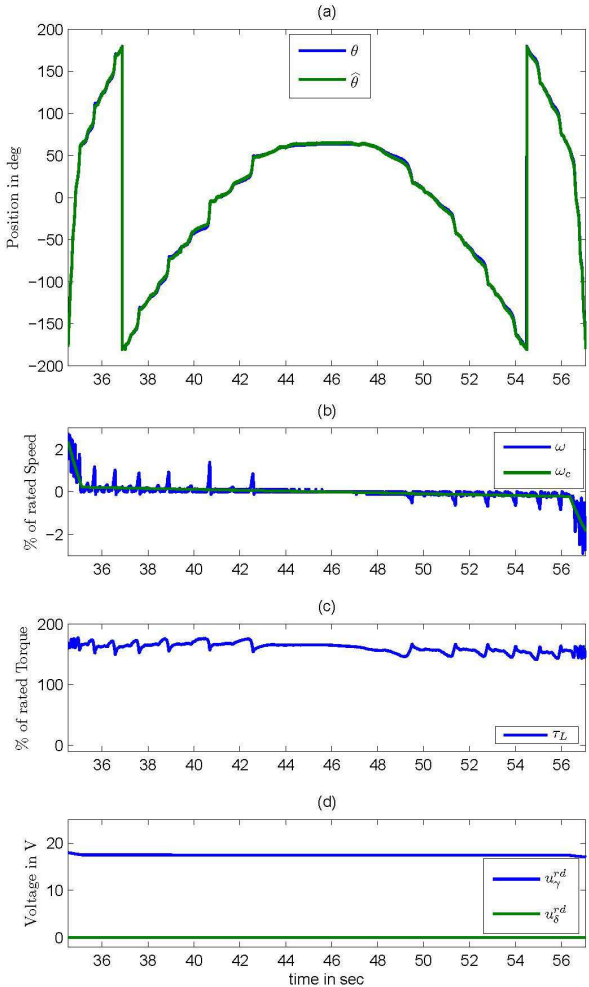}
\caption{Slow speed reversal for IPM:
(a)~measured $\theta$, estimated $\widehat\theta$; (b) measured speed $\omega$, reference speed $\omega_c$; (c) load torque~$\tau_L$; (d) voltages $u_{\gamma\delta}^{rd}$}
\label{fig:SpeedReversalIPM}
\end{figure}
\begin{figure}[b!]
\centering
\includegraphics[width=0.9\columnwidth]{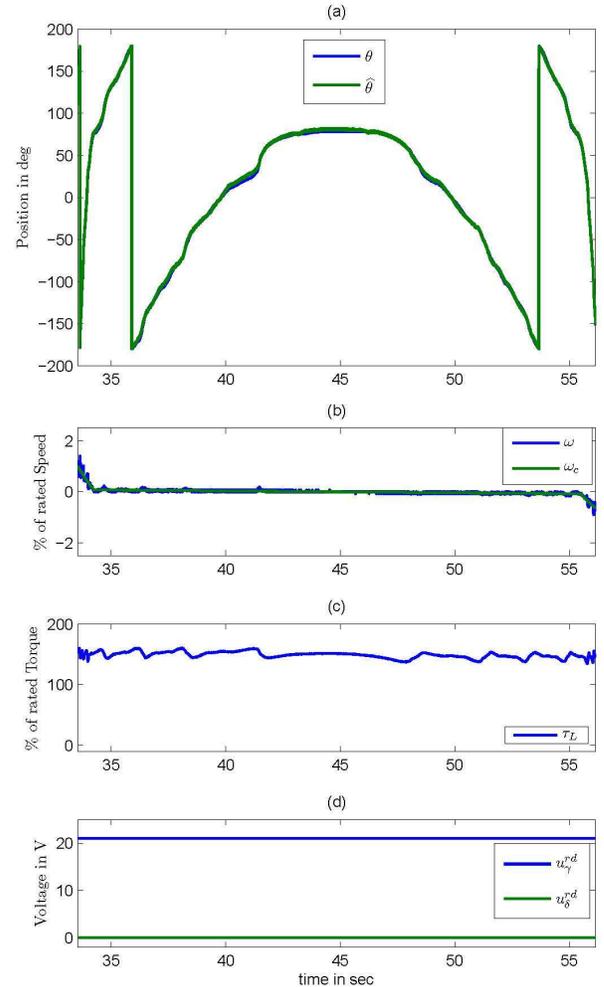}
\caption{Slow speed reversal for SPM:
(a) measured $\theta$, estimated $\widehat\theta$; (b) measured speed $\omega$, reference speed $\omega_c$; (c) load torque~$\tau_L$; (d) voltages $u_{\gamma\delta}^{rd}$}
\label{fig:SpeedReversalSPM}
\end{figure}

\subsection{Validation of the rotor position estimation procedure}
The relevance of the position estimation methodology developed in section~~\ref{sec:pos} is now illustrated on the two test motors, using the parameters estimated in the previous section.
Since the goal is only to test the validity of the angle estimation procedure, a very simple $V/f$ open-loop (i.e. $\Omega_c$ and $\mathcal U_{\gamma\delta}$ do not depend on $i_{\gamma\delta}$) control law is used for~\eqref{eq:claw1}--\eqref{eq:claw5}; a fast-varying ($\Omega:=2\pi\times500\,\text{rad/s}$) square voltage with constant amplitude is added in accordance with~\eqref{eq:HFvoltage}, resulting in
\begin{align*}
    \frac{d\theta_c}{dt} &= \w_c(t)\\
       u_{\gamma\delta} &= u_{\gamma\delta}^{rd}(t) + \w_c(t)\phi_{m} + \widetilde u_{\gamma\delta} f(\Omega t).
\end{align*}
Here $\w_c(t)$ is the motor speed reference; $u_{\gamma\delta}^{rd}(t)$ is a filtered piece-wise constant vector compensating the resistive voltage drop in order to maintain the torque level and the motor stability;
finally $\widetilde u_{\gamma\delta}:=(\widetilde u,0)^T$ with $\widetilde u:=15$\,V.

The rotor position $\widehat\theta$ is then estimated according to~\eqref{eq:nonlin}.

\subsubsection{Long test under various conditions, Fig.~\ref{fig:LongTestIPM}-\ref{fig:LongTestSPM}}
Speed and torque are changed over a period of $210$ seconds; the speed remains between~$\pm5\%$ of the rated speed and the torque varies from $0\%$ to $180\%$ of the rated toque. This represents typical operation conditions at low speed.

When the saturation model is used the agreement between the estimated position $\widehat{\theta}$ and the measured position $\theta$ is very good, with an error always smaller than a few (electrical) degrees. By contrast the estimated error without using the saturation model (i.e. with all the magnetic saturation parameters $\alpha_{ij}$ taken to zero) can reach up to $40^\circ$ for the IPM and $70^\circ$ the SPM. This demonstrates the importance of considering an adequate saturation model including in particular cross-saturation.

\subsubsection{Slow speed reversal, Fig.~\ref{fig:SpeedReversalIPM}-\ref{fig:SpeedReversalSPM}}
This is an excerpt of the long experiment between $35$\,s and~$55$\,s.
The speed is slowly changed from $-0.2\%$ to $+0.2\%$ of the rated speed at~$150\%$ of the rated torque. This is a very demanding test since the motor always remains in the poor observability region, moreover under high load. Once again the estimated angle closely agrees with the measured angle.

\begin{figure}[b!]
\centering
\includegraphics[width=0.9\columnwidth]{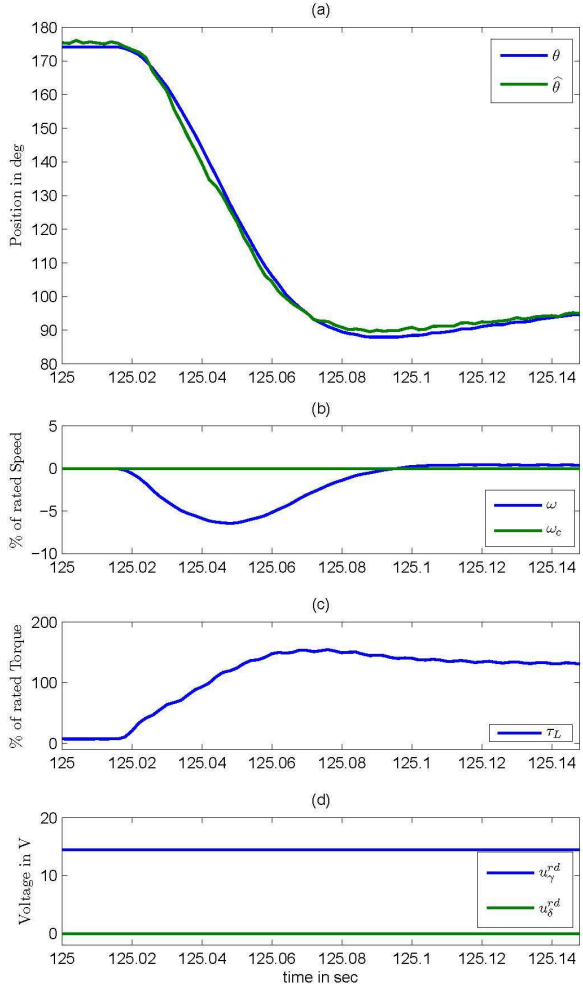}
\caption{Load step at zero speed for IPM:
(a) measured $\theta$, estimated $\widehat\theta$; (b) measured speed $\omega$, reference speed $\omega_c$; (c) load torque~$\tau_L$; (d) voltages $u_{\gamma\delta}^{rd}$}
\label{fig:LoadStep0IPM}
\end{figure}

\begin{figure}[b!]
\centering
\includegraphics[width=0.9\columnwidth]{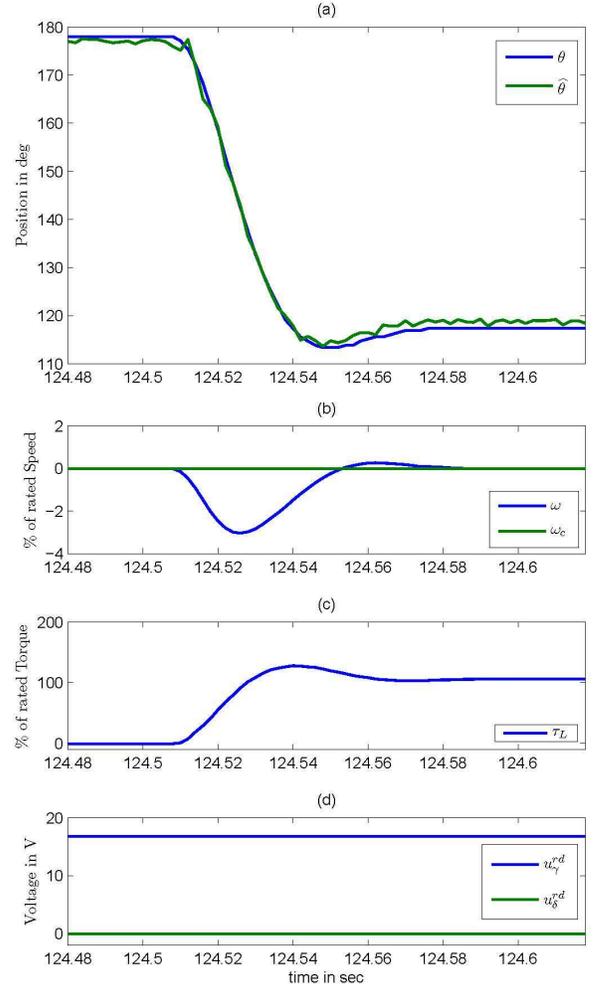}
\caption{Load step at zero speed for SPM:
(a) measured $\theta$, estimated $\widehat\theta$; (b) measured speed $\omega$, reference speed $\omega_c$; (c) load torque~$\tau_L$; (d) voltages $u_{\gamma\delta}^{rd}$}
\label{fig:LoadStep0SPM}
\end{figure}

\subsubsection{Load step at zero speed, Fig.~\ref{fig:LoadStep0IPM}-\ref{fig:LoadStep0SPM}}
This is an excerpt of the long experiment around $t=125$\,s.
The load is suddenly changed from $0\%$ to $100\%$ of the rated torque while the motor is at rest. This test illustrates the quality of the estimation also under dynamic conditions.

\section{Conclusion}
We have presented a simple parametric model of the saturated PMSM together with a new procedure based on signal injection for estimating the rotor angle at low speed relying on an original analysis based on second-order averaging. This is not an easy problem in view of the observability degeneracy at zero speed. The method is general in the sense it can accommodate virtually any control law, saturation model, and form of injected signal. The relevance of the method and the importance of using an adequate magnetic saturation model has been experimentally demonstrated on a SPM motor with little geometric saliency as well as on an IPM motor.

\end{document}